\documentclass[a4paper,10pt]{amsart}
\usepackage[width=0.8\paperwidth,height=0.85\paperheight,twoside=false,includehead,includefoot]{geometry}
\usepackage{amsmath,amssymb,amsfonts}
\usepackage{amsthm,graphicx}


\theoremstyle{plain}
\newtheorem{thm}{Theorem}[section]
\newtheorem{defn}{Definition}

\newtheorem{ex}[thm]{Example}
\newtheorem{rem}[thm]{Remark}

\def\N{\mathbb {N}} \def\Z{\mathbb {Z}} \def\Q{\mathbb {Q}}
\def\R{\mathbb {R}} 
    
\def\d{\partial} 
   
\def\hh{\ast} 

\def\A{\mathcal{A}}
 
\def\k{\textrm{\textup{\textmd{\textbf{k}}}}}  \def\r{\textrm{\textup{\textmd{\textbf{r}}}}}

\def\sg{\sigma}
\def\ve{\varepsilon}   
\def\z{\zeta}

\def\ZZ{\mathcal {Z}}
\def\fh{\frak H}

\def\za{\zeta_{\mathcal{A}}} 
\def\zf{\zeta_{\mathcal{F}}}
\def\zsym{\zeta_{\mathcal{S}}} 

\def\ZA{\mathit{Z}_{\mathcal{A}}}
\def\ZS{\mathit{Z}_{\mathcal{S}}}
\def\ZF{\mathit{Z}_{\mathcal{F}}}
\def\ch{\vee}

\DeclareMathOperator{\wt}{\mathrm{wt}}
\DeclareMathOperator{\dep}{\mathrm{dep}}
\DeclareMathOperator{\h}{\mathrm{ht}}

\begin{document}
\title[On MZVs and FMZVs of maximal height]{On multiple zeta values and finite multiple zeta values of maximal height}
\author{Hideki Murahara and Mika Sakata}
\date{2017.7.6}

\begin{abstract}
An explicit formula for the height-one multiple zeta values was proved by Kaneko and the second author. 
We give an alternative proof of this result and its generalization. 
We also prove its counterpart for the finite multiple zeta values.   
\end{abstract}
\maketitle

\section{Introduction/Main theorem}
For positive integers $k_1, \ldots , k_{d}\in \N$ with $k_1 \geq 2$, the multiple zeta value (MZV) is defined by 
\begin{align*}
\z (k_1,\ldots, k_d) := \sum_{n_1>\cdots >n_d \geq 1} \frac {1}{n_1^{k_1}\cdots n_d^{k_d}}.
\end{align*} 
The quantities $\wt (\k):=k_1+\cdots+k_d$, $\dep (\k):=d$, and $\h (\k):=\#\{ i \,\mid \,k_i\geq 2, 1\le i\leq d\}$ are called the weight, the depth, and the height of the index set $\k =(k_1,\ldots, k_d)$ (or of the multiple zeta value $\z (\k )=\z (k_1,\ldots,k_d)$), respectively. 
For two indices $\k =(k_1,\ldots ,k_i)$ and $\r =(r_1,\ldots ,r_i)$ of the same depth,
$\z (\k +\r )$ denotes $\z (k_1+r_1,\ldots ,k_i+r_i)$. 
In \cite{kaneko_sakata_2016}, M.\ Kaneko and the second author proved the following explicit formula for the height-one multiple zeta values. 
\begin{thm}[\cite{kaneko_sakata_2016}] \label{kaneko-sakata}
For $k,r \in\N$, we have 
\begin{align*}
\z (k+1,\underbrace{1, \ldots , 1}_{r-1})
=\sum^{\min (k,r)}_{i=1}(-1)^{i-1}\sum_{\substack{\wt (\k)=k, \,\wt (\r)=r \\ \dep (\k)=\dep (\r)=i}}
\z (\k+\r). 
\end{align*}
\end{thm}
We note that the right-hand side of this formula is symmetric in $r$ and $k$, and thus the formula makes the duality $\z (k+1,\underbrace{1, \ldots , 1}_{r-1})=\z (r+1,\underbrace{1, \ldots ,1}_{k-1})$ visible.
One of the aims of this paper is to prove a generalization of Theorem \ref{kaneko-sakata}. 
For two indices $\k, \k'$, we say $\k'$ refines $\k$ (denoted  $\k' \succeq \k$) if $\k$ can be obtained from $\k'$ by combining some of its adjacent parts. For example, 
$(1, 2, 3, 4)\succeq (1+2, 3+4)=(3, 7)$. 
Our main theorem is the following: 
\begin{thm} \label{mainthm1}
For $\k=(k_1,\ldots,k_d)\in \N^d$ and $r\in \N$ with $r\geq d$, we have 
\begin{align*} 
\sum_{\substack{r_1+\cdots+r_d=r \\ r_i \geq 1 \, (1\leq i\leq d)}} \z (k_1+1,\underbrace{1, \ldots , 1}_{r_1-1},\ldots,k_d+1,\underbrace{1, \ldots , 1}_{r_d-1})
=\sum_{\substack{\k' \succeq \k\\\dep (\k')\leq r}} \sum_{\substack{\wt (\r)=r \\ \,\dep (\r)=\dep (\k')}} (-1)^{\dep (\k')-d} \z (\k'+\r).  
\end{align*}
\end{thm}
The case $d=1$ of Theorem \ref{mainthm1} gives Theorem \ref{kaneko-sakata}. 
In this paper, we will give two different proofs of Theorem \ref{mainthm1}. 
The first proof is based on the duality theorem and Ohno's relations, and the second is on the derivation relations.   
We also introduce a counterpart of this theorem for the finite multiple zeta values. 
The symbol $\zf$ in the following stands for the finite multiple zeta values, which will be defined in section $4$.  
\begin{thm} \label{mainthm2-1}
For $\k=(k_1,\ldots,k_d)\in \N^d$ and $r\in \N$ with $r\geq d$, we have 
\begin{align*} 
&\sum_{\substack{r_0+\cdots+r_d=r+1 \\ r_i \geq 1 \, (0\leq i\leq d)}} \zf (\underbrace{1, \ldots ,1}_{r_0-1}, k_1+1,\underbrace{1, \ldots , 1}_{r_1-1},\ldots,\underbrace{1, \ldots , 1}_{r_{d-1}-1}, k_d+1,\underbrace{1, \ldots , 1}_{r_d-1}) \\
&=\sum_{\substack{\k' \succeq \k\\\dep (\k')\leq r}} 
\sum_{\substack{\wt (\r)=r \\ \,\dep (\r)=\dep (\k')}} (-1)^{\dep (\k')-d} \zf (\k'+\r). 
\end{align*}
\end{thm}
We note that one must include strings that start with 1's in the case of $\zf$. 
Let us give some examples of the above theorems. 
\begin{ex}
When $k=3,r=4$ in Theorem \ref{kaneko-sakata}, or $\k =(3),r=4$ in Theorem \ref{mainthm1}, we have 
\begin{align*}
\z (4,1,1,1)=&\z (7)-(\z (5,2)+2\z (4,3)+2\z (3,4)+\z (2,5)) \\
&+\z(3,2,2)+\z(2,3,2)+\z(2,2,3). 
\end{align*}  
\end{ex}
\begin{ex}
When $\k =(3,2),r=4$ in Theorem \ref{mainthm1}, we have 
\begin{align*}
&\z (4,3,1,1)+\z (4,1,3,1)+\z (4,1,1,3)\\
&=\z (6,3)+\z (5,4)+\z (4,5) \\
&\,\,\, -(\z (5,2,2)+\z (4,3,2)+2\z (4,2,3)+2\z (3,3,3)+\z (3,2,4)+\z (2,4,3)+\z (2,3,4)) \\
&\,\,\, +\z(3,2,2,2)+\z(2,3,2,2)+\z(2,2,2,3). 
\end{align*}  
\end{ex}
\begin{ex}
When $\k =(3,2),r=4$ in Theorem \ref{mainthm2-1}, we have 
\begin{align*}
&\zf (4,3,1,1)+\zf (4,1,3,1)+\zf (4,1,1,3)+\zf (1,4,3,1)+\zf (1,4,1,3)+\zf (1,1,4,3)\\
&=\zf (6,3)+\zf (5,4)+\zf (4,5) \\
&\,\,\, -(\zf (5,2,2)+\zf (4,3,2)+2\zf (4,2,3)+2\zf (3,3,3)+\zf (3,2,4)+\zf (2,4,3)+\zf (2,3,4)) \\
&\,\,\, +\zf (3,2,2,2)+\zf (2,3,2,2)+\zf (2,2,2,3). 
\end{align*}  
\end{ex}

\begin{rem}
S.~Yamamoto also obtained a generalization of Theorem \ref{kaneko-sakata} by 
using generating function.
\end{rem}

\section{Proof of main theorem}
\subsection{The duality theorem and Ohno's relations for the MZVs} 
We recall some known results that will be used in the proof of Theorem \ref{mainthm1}. 
\begin{defn} \label{du}  
Let $\k=(k_1,\ldots ,k_d)$ be an index set with $k_1 \geq 2$. 
We write 
\[\k=(a_1+1, \underbrace{1,\ldots ,1}_{b_1-1},\ldots ,a_s+1, \underbrace{1,\ldots ,1}_{b_s-1}) \] 
with $a_p, b_q\geq 1$. 
Then, we define the dual index set of $\k$ as 
\[ \k^{\hh}= (b_s+1, \underbrace{1,\ldots ,1}_{a_s-1},\ldots ,b_1+1, \underbrace{1,\ldots ,1}_{a_1-1}). \]
\end{defn}
%
The following result, which is a direct consequence of the iterated integral expression, provides the so-called duality theorem for the MZVs.
\begin{thm} \label{thm:dualitym} 
For $\k=(k_1,\ldots,k_d)$ with $k_1 \geq 2$, we have
\[ \z (\k)=\z (\k^{\hh}). \]
\end{thm}
Ohno's relations is one of the most general explicit relations among the MZVs.
This is a generalization of both the duality theorem above and the well-known sum formula.
\begin{thm}[Ohno \cite{ohno-99}] \label{ohno} 
For $\k =(k_1,\ldots,k_d) \in \N^d$ with $k_1 \geq 2$ and $m\in\Z_{\geq 0}$, we have
\begin{multline*}
\sum_{\substack{\ve_1 + \cdots + \ve_d =m \\ \ve_i \geq 0 \, (1\leq i \leq d)}}
\z (k_1+\ve_1,k_2+\ve_2,\ldots,k_d+\ve_d) 
=\sum_{\substack{\ve_1' + \cdots + \ve_{d'}' = m \\ \ve_i' \geq 0  \, (1\leq i \leq d') }}
\z (k_1'+\ve_1',k_2'+\ve_2',\cdots,k_{d'}'+\ve_{d'}'), 
\end{multline*}
where $(k'_1,\ldots ,k'_{d'})$ is the dual of $\k$. 
\end{thm}

\subsection{Proof of Theorem \ref{mainthm1}} 
We set $\wt (\k):=k$.
\begin{align*}
R.H.S.&=\sum_{\substack{\k' \succeq \k\\\dep (\k')\leq r}} \sum_{\substack{\wt (\r)=r \\ \,\dep (\r)=\dep (\k')}} (-1)^{\dep (\k')-d} \z (\k'+\r) \\
&=\sum_{j=d}^{\min(k, r)}(-1)^{j-d}\sum_{(k'_{1}, \ldots , k'_{j}) \succeq \k}
\sum_{\substack{r_{1}+\cdots +r_j=r\\ r_i \geq 1 \, (1\leq i\leq j) }} \z (k'_{1}+r_{1}, \ldots, k'_{j}+r_{j}).
\end{align*}
From Ohno's relations, we have  
\begin{align*}
&\sum_{\substack{r_{1}+\cdots +r_{j}=r \\ r_i \geq 1 \, (1\leq i\leq j) }} \z (k'_{1}+r_{1}, \ldots, k'_{j}+r_{j}) \\
&=\sum_{\substack{r_{1}+\cdots +r_{j}=r-j \\ r_i \geq 0 \, (1\leq i\leq j) }} \z (k'_{1}+r_{1}+1, \ldots, k'_{j}+r_{j}+1) \\
&=\sum_{\substack{r_{1}+\cdots +r_{k}=r-j \\ r_i \geq 0 \, (1\leq i\leq k) }} \z (r_{1}+2, \underbrace{r_{2}+1, \ldots ,  r_{k'_j}+1}_{k'_{j}-1},\ldots, r_{k-k'_{1}+1}+2, \underbrace{r_{k-k'_{1}+2}+1, \ldots, r_{k}+1}_{k'_{1}-1}). 
\end{align*}
Then, 
\begin{align} 
\begin{split} \label{g}
R.H.S.&=\sum_{j=d}^{\min(k, r)}(-1)^{j-d}\sum_{(k'_{1}, \ldots , k'_{j}) \succeq \k}
 \sum_{\substack{r_{1}+\cdots +r_{k}=r-j \\ r_i \geq 0 \, (1\leq i\leq k) }}  
 \z (r_{1}+2, \underbrace{r_{2}+1, \ldots r_{k'_j}+1}_{k'_{j}-1}, \ldots  \\
&\qquad\qquad\qquad\qquad\qquad\qquad\qquad\qquad\quad\,\,\,  \ldots , r_{k-k'_{1}+1}+2, \underbrace{r_{k-k'_{1}+2}+1, \ldots, r_{k}+1}_{k'_{1}-1}). 
\end{split}
\end{align}
Now, for a fixed $\k=(k_1,\ldots ,k_d) \in \N^d$ and $r\in\N$ with $r\geq d$, and a variable $s\in\N$, we set 
\begin{align*}
F(s)&:=F_{\k,r}(s) \\
&:=\sum_{\substack{r_{1}+\cdots+r_{k}=r-d\\ r_{i}\geq 0\, (1\leq i \leq k), \, \h =s}}\zeta(r_{1}+2,\underbrace{r_{2}+1,\ldots,r_{k_{d}}+1}_{k_{d}-1}, \ldots, r_{k-k_{1}+1}+2,\underbrace{r_{k-k_{1}+2}+1,\ldots, r_{k}+1}_{k_{1}-1}),  
\end{align*}
where the sum runs over those $r_{i}$s such that the argument of $\z$ is of height $s$. 
For each $j$ in equation (\ref{g}), the number of each height $s$ MZV is $\binom{s-d}{j-d}$ because it is determined by the partitions of $\k$ and the allocation of $r$. 
(Here, we note that there are $s-d$ places with a component greater than 1 that are determined by the partitions of $\k$ and the allocation of $r$, and there are $j-d$ places which depend on the various partitions of $\k$.) 
Thus, by focusing on the height and retake the sums, we have
\begin{align*} 
R.H.S.=\sum_{j=d}^{\min(k,r)} (-1)^{j-d} \sum_{s=j}^{\min(k,r)} \binom{s-d}{j-d} F(s). 
\end{align*}
From the binomial theorem and the duality theorem, we have 
\begin{align*}
R.H.S.&=\scalebox{0.89}[1]{$\displaystyle   
 \binom{0}{0} F(d)+\left( \binom{1}{0}-\binom{1}{1} \right) F(d+1)+\cdots +\left( \sum_{m=0}^{\min(k,r)-d}
 (-1)^m \binom{\min(k,r)-d}{m} \right) F(\min(k,r)) 
 $} \\
&=F(d) \\
&=\sum_{\substack{r_{1}+\cdots+r_{d}=r-d \\ r_{i}\geq 0 \, (1\leq i\leq d) }}
\zeta(r_{1}+2,\underbrace{1,\ldots,1}_{k_{d}-1},\ldots,r_{d}+2,\underbrace{1,\ldots,1}_{k_{1}-1}) \\ 
&=\sum_{\substack{r_{1}+\cdots+r_{d}=r \\ r_{i}\geq 1 \, (1\leq i\leq d) }}
\zeta(r_{1}+1,\underbrace{1,\ldots,1}_{k_{d}-1},\ldots, r_{d}+1,\underbrace{1,\ldots,1}_{k_{1}-1}) \\
&=\sum_{\substack{r_1+\cdots+r_d=r \\ r_i \geq 1 \, (1\leq i\leq d)}} \z (k_1+1,\underbrace{1, \ldots , 1}_{r_1-1},\ldots,k_d+1,\underbrace{1, \ldots , 1}_{r_d-1}) \\
&=L.H.S. 
\end{align*}

\section{Alternative proof}
\subsection{Algebraic setup and the derivation relations}
In this section, we give an alternative proof of Theorem \ref{mainthm1}. 
Throughout the proof, we use the algebraic setup introduced by M.\ Hoffman in \cite{hoffman-97}. 
Let $\fh :=\Q \left\langle x,y \right\rangle$ be the noncommutative polynomial ring in two indeterminates $x$, $y$, and $\fh^1$ (resp. $\fh^0$) its subrings $\Q +\fh y$ (resp. $\Q +x \fh y$). 
We set $z_{k} := x^{k-1} y$ $(k\in\N)$. 
Then $\fh^1$ is freely generated by $\{z_{k}\}_{k \geq 1}$. 
We define the $\Q$-linear map $\mathit{Z}:\fh^0 \to \R$ by $\mathit{Z}(1):=1$ and $\mathit{Z}( z_{k_1} \cdots z_{k_d}):= \z (k_1,\ldots, k_d)$. 

A derivation on $\fh$ is a $\Q$-linear map $\partial:\frak{H}\rightarrow\frak{H}$ satisfying Leibniz's rule $\partial(ww^{\prime})=\partial(w)w^{\prime}+w\partial(w^{\prime})$. Such a derivation is uniquely determined by its images of generators $x$ and $y$. Set $z:=x+y$. 
For each $l\in\N$, the derivation $\d_l$ on $\fh$ is defined by $\partial_l(x):=xz^{l-1}y$ and $\d_l(y):=-xz^{l-1}y$. 
We note that $\d_l(1)=0$ and $\d_l(z)=0$. 
K.~Ihara, Kaneko and D.~Zagier proved the derivation relations for the MZVs.
\begin{thm}[Ihara--Kaneko--Zagier \cite{ihara_kaneko_zagier_2006}] \label{ihara-kaneko-zagier}
For $l\in\N$, we have 
\begin{align*} 
 \mathit{Z} (\partial_{l}(w)) = 0 \quad (w \in \fh^0).  
\end{align*}
\end{thm}

Now, we prepare some additional notation. 
Let $\alpha$ be the endmorphism on $\fh$ such that $\alpha (x):=x-xy$ and $\alpha (y):=-xy$, 
and $\tau$ be the anti-automorphism on $\fh$ such that $\tau (x):=y$ and $\tau (y):=x$. 
For positive integer $m$, we also let $\beta_m:\frak{H}\rightarrow\frak{H}$ by setting $\beta_{m}(w)$ to be the degree $m$ part of $w\in\frak{H}$. 
For  each $l\in\N$, we define the derivation $D_l$ on $\fh$ by $D_l(x):=0$ and $D_l(y):=x^ly$. 
Set $\sg:=\textrm{exp} \left( \sum_{l=1}^{\infty}\frac{D_l}{l} \right)$. 
Then, we find the map $\sg$ is an automorphism on $\widehat{\fh}=\Q\left\langle\left\langle x,y\right\rangle\right\rangle$, and we see that $\sg(x)=x$ and $\sg(y)=\frac{1}{1-x}y$. 
(See also \cite[Section $6$]{ihara_kaneko_zagier_2006} and \cite[Appendix]{tanaka_2009}.)

\subsection{Alternative proof of Theorem \ref{mainthm1}}
According to Ihara, Kaneko and Zagier in \cite[Proof of Theorem 3]{ihara_kaneko_zagier_2006}, we first note that  
\begin{align} \label{5}
\sg-\tau\sg\tau=\left( 1-\textrm{exp} \, \left( \sum_{l=1}^{\infty} \frac{\d_l}{l}\right) \right) \sg . 
\end{align}
By Theorem \ref{ihara-kaneko-zagier}, we have
\[
 \mathit{Z} \left( \beta_m \left( \left( 1-\textrm{exp} \, \left( \sum_{l=1}^{\infty} \frac{\d_l}{l} \right) \right) \sg (w) \right) \right) =0 \quad (m\in\mathbb{N}, w\in\fh^0 ). 
\]
Let $k:=k_1+\cdots +k_d$.
From the equality (\ref{5}) and by putting $\alpha (x^{k_1-1}y\cdots x^{k_d-1}y)$ into $w$,
\[ \mathit{Z} (\beta_m (\sg-\tau\sg\tau)\alpha (x^{k_1-1}y\cdots x^{k_d-1}y))=0 \quad (m\ge k). \]
Here, 
\begin{align*} 
\beta_{k+r} \sigma\alpha (x^{k_1-1}y\cdots x^{k_d-1}y)
&=\beta_{k+r} \sigma ((x-xy)^{k_1-1}(-xy)\cdots (x-xy)^{k_d-1}(-xy)) \\
&=(-1)^d \beta_{k+r} \sigma ((x-xy)^{k_1-1}xy \cdots (x-xy)^{k_d-1}xy) \\ 
&=(-1)^d \beta_{k+r} \biggl( \left( x-\frac{x}{1-x}y\right)^{k_1-1}\frac{x}{1-x}y \cdots \left( x-\frac{x}{1-x}y\right)^{k_d-1}\frac{x}{1-x}y \biggr)\\
&=(-1)^d \sum_{ \substack{e_{1,1}+\cdots+e_{d,k_d}=r \\ e_{i,j} \geq 0 \, (1\leq j\leq k_i-1) \\ e_{i,j} \geq 1 \, (j=k_i)} }
 (-x^{e_{1,1}}y)\cdots (-x^{e_{1,{k_1-1}}}y) x^{e_{1,k_1}}y \\ 
&\qquad\qquad\qquad\qquad\quad\,\,\, \cdot\cdots \cdots\cdot (-x^{e_{d,1}}y)\cdots (-x^{e_{d,k_{d}-1}}y) x^{e_{d,k_d}}y. 
\end{align*}
When  $e_{i,j}=0$, we understand $x^{e_{i,j}}y=-x$. 
On the other hand, 
\begin{align*}
\beta_{k+r} \tau\sigma\tau\alpha (x^{k_1-1}y\cdots x^{k_d-1}y)
&=\beta_{k+r} \tau\sigma\tau ((x-xy)^{k_1-1}(-xy)\cdots (x-xy)^{k_d-1}(-xy)) \\ 
&=(-1)^d \beta_{k+r} \tau\sigma\tau ((x-xy)^{k_1-1}xy \cdots (x-xy)^{k_d-1}xy) \\ 
&=(-1)^d \beta_{k+r} \tau\sigma (xy (y-xy)^{k_d-1}\cdots xy(y-xy)^{k_1-1}) \\ 
&=(-1)^d \beta_{k+r} \tau \left( \frac{x}{1-x}y^{k_d} \cdots \frac{x}{1-x}y^{k_1} \right) \\ 
&=(-1)^d \beta_{k+r} \left( x^{k_1}\frac{y}{1-y} \cdots x^{k_d}\frac{y}{1-y} \right) \\
&=(-1)^d \sum_{\substack{r_1+\cdots+r_d=r \\ r_i \geq 1 \, (1\leq i\leq d)}} x^{k_1}y^{r_1}\cdots x^{k_d}y^{r_d}. 
\end{align*}
Then, we have 
\begin{align*}
&\sum_{ \substack{e_{1,1}+\cdots+e_{d,k_d}=r \\ e_{i,j} \geq 0 \, (1\leq j\leq k_i-1) \\ e_{i,j} \geq 1 \, (j=k_i)} }
 \mathit{Z} ( (-x^{e_{1,1}}y)\cdots (-x^{e_{1,{k_1-1}}}y) x^{e_{1,k_1}}y \\ 
&\qquad\qquad\qquad\quad\,\,\, \cdot\cdots \cdots\cdot (-x^{e_{d,1}}y)\cdots (-x^{e_{d,k_{d}-1}}y) x^{e_{d,k_d}}y ) \\
&=\sum_{\substack{r_1+\cdots+r_d=r \\ r_i \geq 1 \, (1\leq i\leq d)}}
 \mathit{Z} ( x^{k_1}y^{r_1}\cdots x^{k_d}y^{r_d} ). 
\end{align*}
In the L.H.S. of the above equality, each $(-x^{e_{1,1}}y)\cdots (-x^{e_{1,{k_i-1}}}y) x^{e_{1,k_i}}y$ is corresponding to `the refinement of $k_i$' $+$ `the part of $\textrm{\textup{\textmd{\textbf{r}}}}$' with signs of the R.H.S. of Theorem \ref{mainthm1}. 
Thus, we have the desired result.

\section{The counterpart of Theorem \ref{mainthm1} for the finite multiple zeta values}
\subsection{Definitions and second main result}
In this section, we prove the counterpart of Theorem \ref{mainthm1} for what we call `finite multiple zeta values (FMZVs)', a generic term for the $\A$-finite multiple zeta values and the symmetrized multiple zeta values.

We consider the collection of truncated sums $\z_{p}(k_1, \ldots, k_d):=\sum_{p>n_1>\cdots>n_d\geq1}\frac{1}{n_1^{k_1}\cdots n_d^{k_d}}$ modulo all primes $p$ in the quotient ring $\A=(\prod_{p}\Z/p\Z)/(\bigoplus_{p}\Z/p\Z)$, which is a $\Q$-algebra. 
Elements of $\A$ are represented by $(a_p)_p$, where $a_p\in\Z/p\Z$, and two elements $(a_p)_p$ and $(b_p)_p$ are identified if and only if $a_p=b_p$ for all but finitely many primes $p$.
For integers $k_1, \ldots ,k_{d} \in \N$, the $\A$-finite multiple zeta value $\za(k_1,\cdots, k_d)$ is defined by
\begin{align*}
 \za (k_1, \ldots, k_d)&:=\biggl( \sum_{p>n_1>\cdots >n_d\geq1}\frac{1}{n_1^{k_1}\cdots n_d^{k_d}} \bmod{p} \biggr)_{p} \in \A.
\end{align*}

The symmetrized multiple zeta values was first introduced by Kaneko and Zagier in \cite{kaneko_2014, kaneko_zagier-2015}. 
For integers $k_1,\ldots ,k_{d} \in \N$, we let  
\begin{align*}
\zsym ^{\ast} (k_1, \ldots, k_d) &:= \sum_{i=0}^{d}(-1)^{k_1+\cdots +k_i}\z^{\ast}(k_i, \ldots , k_1) \z ^{\ast}(k_{i+1}, \ldots , k_d). 
\end{align*}
Here, the symbols $\z ^{\ast}$ on the right-hand sides stand for the regularized values coming from harmonic regularizations,
i.e.,  real values obtained by taking constant terms of harmonic regularizations as explained in   \cite{ihara_kaneko_zagier_2006}. 
In the sum, we understand $\zeta^{\ast}(\emptyset )=1$.
Let $\ZZ_{\R}$ be the $\Q$-vector subspace of $\R$ spanned by $1$ and all MZVs, which is 
a $\Q$-algebra. 
Then, the symmetrized multiple zeta value $\zsym (k_1, \ldots , k_d)$ is defined as an element 
in the quotient ring $\ZZ_{\R}/\zeta (2)$ by 
\[  \zsym (k_1, \ldots , k_d) := \zsym^{\ast} (k_1, \ldots , k_d) \bmod \zeta(2). \]
(For more details, see \cite{kaneko_2014, kaneko_zagier-2015}.)

Now, we restate our second main result which was introduced in section 1. 
\begin{thm} \label{mainthm2}
For $\k=(k_1,\ldots,k_d)\in \N^d$ and $r\in \N$ with $r\geq d$, we have 
\begin{align*} 
&\sum_{\substack{r_0+\cdots+r_d=r+1 \\ r_i \geq 1 \, (0\leq i\leq d)}} \zf (\underbrace{1, \ldots ,1}_{r_0-1}, k_1+1,\underbrace{1, \ldots , 1}_{r_1-1},\ldots,\underbrace{1, \ldots , 1}_{r_{d-1}-1}, k_d+1,\underbrace{1, \ldots , 1}_{r_d-1}) \\
&=\sum_{\substack{\k' \succeq \k\\\dep (\k')\leq r}} 
\sum_{\substack{\wt (\r)=r \\ \,\dep (\r)=\dep (\k')}} (-1)^{\dep (\k')-d} \zf (\k'+\r) 
\qquad\quad (\mathcal{F}=\mathcal{A} \textrm{ or } \mathcal{S}). 
\end{align*}
\end{thm}

\begin{rem}
Denoting $\ZZ_{\A}$ by the $\Q$-vector subspace of $\A$ spanned by $1$ and all $\A$-finite multiple zeta values, Kaneko and Zagier conjecture that the homomorphism $\ZZ_{\A} \rightarrow \ZZ_{\R}/\z(2)$ sending $\za (k_1,\ldots,k_d)$ to $\zsym (k_1,\ldots,k_d)$ is an isomorphism.   
\end{rem}

\subsection{Proof of Theorem \ref{mainthm2}}
\subsubsection{Ohno type relations for the FMZVs}
K.\ Oyama proved Ohno type relations for the FMZVs, which was first conjectured by Kaneko in \cite{kaneko_2014}. 
\begin{defn}
For $\k =(k_1,\ldots,k_d) \in \N^d$, we define Hoffman's dual of $\k$ by
\begin{equation*}
\k^{\ch}:=(\underbrace{1,\ldots,1}_{k_1}+\underbrace{1,\ldots,1}_{k_2}+1,\ldots,1+\underbrace{1,\ldots,1}_{k_d}).
\end{equation*}
\end{defn}

Then, Ohno type relations is the following:
\begin{thm}[Oyama \cite{oyama_2015}]
For $(k_1,\ldots,k_d) \in \N^d$ and $m \in \Z_{\geq 0}$, we have
\begin{equation*}
\sum_{\substack{\ve_1 + \cdots + \ve_d =m \\ \ve_i \geq 0 \, (1\leq i \leq d)}} 
\zf (k_1+\ve_1,\ldots,k_d+\ve_d)
=\sum_{\substack{\ve_1' + \cdots + \ve_{d'}' = m \\ \ve_i' \geq 0 \, (1\leq i \leq d') }}
\zf ((k'_1+\ve'_1,\ldots,k'_{d'}+\ve'_{d'})^{\ch}),
\end{equation*}
where $(k'_1,\ldots,k'_{d'})=(k_1,\ldots,k_d)^{\ch}$ is Hoffman's dual of $(k_1,\ldots,k_d)$.
\end{thm}

\subsubsection{Proof of Theorem \ref{mainthm2}}
We can prove Theorem \ref{mainthm2} in the same manner as in the proof of Theorem \ref{mainthm1}. 
From Ohno type relations, we have
\begin{align*}
&R.H.S.=\sum_{\substack{\k' \succeq \k\\\dep (\k')\leq r}} \sum_{\substack{\wt (\r)=r \\ \,\dep (\r)=\dep (\k')}} (-1)^{\dep (\k')-d} \zf (\k'+\r) \\
&=\sum_{j=d}^{\min(k, r)}(-1)^{j-d}\sum_{(k'_{1}, \ldots , k'_{j}) \succeq \k}
\sum_{\substack{r_{1}+\cdots +r_j=r-j \\ r_i \geq 0 \, (1\leq i\leq j) }} \zf (k'_{1}+r_{1}+1, \ldots, k'_{j}+r_{j}+1) \\
&=\sum_{j=d}^{\min(k, r)}(-1)^{j-d}\sum_{(k'_{1}, \ldots , k'_{j}) \succeq \k}
 \sum_{\substack{r_0+\cdots +r_k=r-j \\ r_i \geq 0 \, (0\leq i\leq k) }}  
 \zf ((\underbrace{r_0+1, \ldots, r_{k'_1-1}+1}_{k'_1}, r_{k'_1}+2,    \\
&\quad\quad\qquad\qquad\qquad\qquad\qquad\,\,\,\, 
 \underbrace{r_{k'_1+1}+1, \ldots, r_{k'_1+k'_2-1}+1}_{k'_2-1}, \ldots, \underbrace{r_{k-k'_{j-1}-k'_j+1}+1, \ldots, r_{k-k'_j-1}+1}_{k'_{j-1}-1},  \\
&\qquad\qquad\qquad\qquad\qquad\qquad\qquad\qquad\qquad\qquad\qquad\qquad\,\,\,\,\,\,
 r_{k-k'_j}+2, \underbrace{r_{k-k'_j+1}+1, \ldots, r_{k}+1}_{k'_j})^{\vee} ).  
\end{align*}
For a fixed $\k=(k_1,\ldots ,k_d) \in \N^d$ and $r\in\N$ with $r \geq d$, and a variable $s\in\N$, we set 
\begin{align*}
G(s):=G_{\k,r}(s) &:= \sum_{\substack{r_0+\cdots +r_{k}=r-d-2 \\ r_i \geq 0 \, (0\leq i\leq k),  \, \h =s}}  
 \zf ((r_0+2, \underbrace{r_1+1, \ldots, r_{k_1-1}+1}_{k_1-1}, r_{k_1}+2,    \\
&\quad\quad\quad\quad\,\,\,\,\,\,\,\, 
 \underbrace{r_{k_1+1}+1, \ldots, r_{k_1+k_2-1}+1}_{k_2-1}, \ldots, \underbrace{r_{k-k_{d-1}-k_d+1}+1, \ldots, r_{k-k_d-1}+1}_{k_{d-1}-1},  \\
&\quad\qquad\qquad\qquad\qquad\qquad\qquad\quad
 r_{k-k_d}+2, \underbrace{r_{k-k_d+1}+1, \ldots, r_{k-1}+1}_{k_d-1}, r_{k}+2)^{\vee} ),  
\end{align*}
\begin{align*}
G'(s):=G'_{\k,r}(s) &:= \sum_{\substack{r_1+\cdots +r_{k}=r-d-1 \\ r_i \geq 0 \, (1\leq i\leq k),  \, \h =s}}  
 \zf ((1, \underbrace{r_1+1, \ldots, r_{k_1-1}+1}_{k_1-1}, r_{k_1}+2,    \\
&\quad\quad\quad\quad\,\,\,\,\,\,\,\, 
 \underbrace{r_{k_1+1}+1, \ldots, r_{k_1+k_2-1}+1}_{k_2-1}, \ldots, \underbrace{r_{k-k_{d-1}-k_d+1}+1, \ldots, r_{k-k_d-1}+1}_{k_{d-1}-1},  \\
&\qquad\qquad\qquad\qquad\qquad\qquad\quad\quad
 r_{k-k_d}+2, \underbrace{r_{k-k_d+1}+1, \ldots, r_{k-1}+1}_{k_d-1}, r_{k}+2)^{\vee} )  \\
&\,\,\,\, +\sum_{\substack{r_0+\cdots +r_{k-1}=r-d-1 \\ r_i \geq 0 \, (0\leq i\leq k-1),  \, \h =s}}  
 \zf ((r_0+2,\underbrace{r_1+1, \ldots, r_{k_1-1}+1}_{k_1-1}, r_{k_1}+2,    \\
&\quad\quad\quad\quad\,\,\,\,\,\,\,\, 
 \underbrace{r_{k_1+1}+1, \ldots, r_{k_1+k_2-1}+1}_{k_2-1}, \ldots, \underbrace{r_{k-k_{d-1}-k_d+1}+1, \ldots, r_{k-k_d-1}+1}_{k_{d-1}-1},  \\
&\quad\qquad\qquad\qquad\qquad\qquad\qquad\qquad\,\,\,\,\,
 r_{k-k_d}+2, \underbrace{r_{k-k_d+1}+1, \ldots, r_{k-1}+1}_{k_d-1}, 1)^{\vee} ),  
\end{align*}
\begin{align*}
G''(s):=G''_{\k,r}(s) &:= \sum_{\substack{r_1+\cdots +r_{k-1}=r-d \\ r_i \geq 0 \, (1\leq i\leq k-1),  \, \h =s}}  
 \zf ((1, \underbrace{r_1+1, \ldots, r_{k_1-1}+1}_{k_1-1}, r_{k_1}+2,    \\
&\quad\quad\qquad\,\,\,\,\,
 \underbrace{r_{k_1+1}+1, \ldots, r_{k_1+k_2-1}+1}_{k_2-1}, \ldots, \underbrace{r_{k-k_{d-1}-k_d+1}+1, \ldots, r_{k-k_d-1}+1}_{k_{d-1}-1},  \\
&\qquad\qquad\qquad\qquad\qquad\qquad\qquad\quad\,\,
 r_{k-k_d}+2, \underbrace{r_{k-k_d+1}+1, \ldots, r_{k-1}+1}_{k_d-1}, 1)^{\vee} ),   
\end{align*}
where the sum runs over those $r_{i}$s such that the argument of $\zf$ is of height $s$. 
By concentrating on the height, and adding up all the terms, then re-arranging the sums; we have 
\begin{align*} 
R.H.S.&=\sum_{j=d+1}^{\min(k+1,r-1)} (-1)^{j-d-1} \sum_{s=j}^{\min(k+1,r-1)} \binom{s-d-1}{j-d-1} G(s) \\ 
&\,\,\,\, +\sum_{j=d}^{\min(k,r-1)} (-1)^{j-d} \sum_{s=j}^{\min(k,r-1)} \binom{s-d}{j-d} G'(s) \\ 
&\,\,\,\, +\sum_{j=d-1}^{\min(k-1,r-1)} (-1)^{j-d+1} \sum_{s=j}^{\min(k-1,r-1)} \binom{s-d+1}{j-d+1} G''(s) \\ 
&=G(d+1)+G'(d)+G''(d-1) \\
&=\sum_{\substack{r_0+\cdots+r_{d}=r-d \\ r_{i}\geq 0 \, (0\leq i\leq d) }}
 \zf ((r_0+1,\underbrace{1,\ldots,1}_{k_{1}-1},r_1+2, \ldots,r_{d-1}+2,\underbrace{1,\ldots,1}_{k_d-1},r_d+1)^{\vee} ) \\ 
&=\sum_{\substack{r_0+\cdots+r_d=r+1 \\ r_i \geq 1 \, (0\leq i\leq d)}} \zf (\underbrace{1, \ldots ,1}_{r_0-1}, k_1+1,\underbrace{1, \ldots , 1}_{r_1-1},\ldots,\underbrace{1, \ldots , 1}_{r_{d-1}-1}, k_d+1,\underbrace{1, \ldots , 1}_{r_d-1}) \\
&=L.H.S. 
\end{align*}

\subsection{Alternative proof of Theorem \ref{mainthm2}}
\subsubsection{The derivation relations for the FMZVs}
The derivation relations for the FMZVs is conjectured by Oyama and proved by the first author in  \cite{murahara_2016}. 

We define two $\Q$-linear maps $\ZA \colon\fh^1 \to \A$ and $\ZS \colon\fh^1 \to \ZZ_{\R}/\zeta (2)$ respectively by $\ZA (1):=1$ and $\ZA (z_{k_1} \cdots z_{k_d}):=\za (k_1,\ldots ,k_d)$, and $\ZS (1):=1$ and $\ZS (z_{k_1} \cdots z_{k_d}):=\zeta_{\mathcal{S}} (k_1,\ldots ,k_d)$. 
We also define the $\Q$-linear operator $L_x$ on $\fh$ by $L_x:=xw \,\, (w\in\fh )$. 
\begin{thm} \label{derf}
For $l\in\N$, we have 
\begin{align*} \label{1}
 \ZF (L_x^{-1}\partial_{l}L_x(w)) = 0 \quad (w \in \fh^1 ,\mathcal{F}=\mathcal{A} \textrm{ or } \mathcal{S}).  
\end{align*}
\end{thm} 
\subsubsection{Alternative proof of Theorem \ref{mainthm2}}
By Theorem \ref{derf}, we have
\[
 \ZF \left( \beta_m \left( L_x^{-1} \left( 1-\textrm{exp} \, \left( \sum_{l=1}^{\infty} \frac{\d_l}{l} \right) \right) L_x \sg (w) \right) \right) =0 \quad (m\in\mathbb{N}, w\in \fh^1 ). 
\]
Since $L_x\sg =\sg L_x$, 
\[
 \ZF \left( \beta_m \left( L_x^{-1} \left( 1-\textrm{exp} \, \left( \sum_{l=1}^{\infty} \frac{\d_l}{l} \right) \right) \sg L_x (w) \right) \right) =0 \quad (m\in\mathbb{N}, w\in \fh^1 ). 
 \]
Let $k:=k_1+\cdots +k_d$.
From the equality (\ref{5}) and by putting $\alpha (x^{k_1-1}y\cdots x^{k_d-1}y)$ into $w$, we have
\[ \mathit{Z}_{\mathcal{F}} (\beta_m L_x^{-1}(\sigma-\tau\sigma\tau)L_x\alpha (x^{k_1-1}y\cdots x^{k_d-1}y))=0 \quad (m\ge k). \]
Since $L_x^{-1}\sg L_x=\sg$ and by the same calculation in subsection 3.2, we have
\begin{align*} 
 \beta_{k+r} L_x^{-1}\sigma L_x\alpha (x^{k_1-1}y\cdots x^{k_d-1}y)
 &=\beta_{k+r} \sigma \alpha (x^{k_1-1}y\cdots x^{k_d-1}y) \\
 &=(-1)^d \sum_{ \substack{e_{1,1}+\cdots+e_{d,k_d}=r \\ e_{i,j} \geq 0 \, (1\leq j\leq k_i-1) \\ e_{i,j} \geq 1 \, (j=k_i)} }
  (-x^{e_{1,1}}y)\cdots (-x^{e_{1,{k_1-1}}}y) x^{e_{1,k_1}}y \\ 
 &\qquad\qquad\qquad\qquad\qquad\quad\,\,\, \cdot\cdots \cdots\cdot (-x^{e_{d,1}}y)\cdots (-x^{e_{d,k_{d}-1}}y) x^{e_{d,k_d}}y.  
\end{align*}
When  $e_{i,j}=0$, we understand $x^{e_{i,j}}y=-x$.  
On the other hand, 
\begin{align*}
 \beta_{k+r} L_x^{-1}\tau\sigma\tau L_x\alpha (x^{k_1-1}y\cdots x^{k_d-1}y)
 &=\beta_{k+r} L_x^{-1}\tau\sigma\tau L_x ((x-xy)^{k_1-1}(-xy)\cdots (x-xy)^{k_d-1}(-xy)) \\ 
 &=(-1)^d \beta_{k+r} L_x^{-1}\tau\sigma\tau L_x ((x-xy)^{k_1-1}xy \cdots (x-xy)^{k_d-1}xy) \\
 &=(-1)^d \beta_{k+r} L_x^{-1}\tau \sigma (xy (y-xy)^{k_d-1}\cdots xy(y-xy)^{k_1-1}y) \\ 
 &=(-1)^d \beta_{k+r} L_x^{-1}\tau \left( \frac{x}{1-x}y^{k_d} \cdots \frac{x}{1-x}y^{k_1} \frac{1}{1-x}y \right) \\ 
 &=(-1)^d \beta_{k+r} \left( \frac{1}{1-y} x^{k_1}\frac{y}{1-y} \cdots x^{k_d}\frac{y}{1-y} \right) \\
 &=(-1)^d \sum_{\substack{r_0+\cdots+r_d=r+1 \\ r_i \geq 1 \, (0\leq i\leq d)}} y^{r_0-1} x^{k_1}y^{r_1}\cdots x^{k_d}y^{r_d}. 
\end{align*}
Then, we have 
\begin{align*}
&\sum_{ \substack{e_{1,1}+\cdots+e_{d,k_d}=r \\ e_{i,j} \geq 0 \, (1\leq j\leq k_i-1) \\ e_{i,j} \geq 1 \, (j=k_i)} }
 \mathit{Z}_{\mathcal{F}} ((-x^{e_{1,1}}y)\cdots (-x^{e_{1,{k_1-1}}}y) x^{e_{1,k_1}}y \\ 
&\qquad\qquad\qquad\qquad \cdot\cdots \cdots\cdot (-x^{e_{d,1}}y)\cdots (-x^{e_{d,k_{d}-1}}y) x^{e_{d,k_d}}y ) \\
&=\sum_{\substack{r_0+\cdots+r_d=r+1 \\ r_i \geq 1 \, (0\leq i\leq d)}} 
 \mathit{Z}_{\mathcal{F}} (y^{r_0-1} x^{k_1}y^{r_1}\cdots x^{k_d}y^{r_d}). 
\end{align*}
This completes the proof of the theorem. 

\section*{Acknowledgment}
The authors would like to express their sincere gratitude to Professor Masanobu Kaneko for valuable comments and advice. The authors also wish to thank Professor Yasuo Ohno and Professor Shingo Saito for their helpful advice. 
The second author is supported by Japan Society for the Promotion of Science, Grant-in-Aid for JSPS Fellows 14J00005 (M.S.).

\bigskip

\begin{flushleft}
\begin{small}
{H.~Murahara}: 
{Nakamura Gakuen University Graduate School, 5-7-1, Befu, Jonan-ku Fukuoka-shi, Fukuoka, 814-0104, Japan}, 
e-mail: {\tt hmurahara@nakamura-u.ac.jp}

\

{M.~Sakata}: 
{Midori Seiho High School, 6-3-9 Ikeshima-cho, Higashiosaka-shi, Osaka 579-8064, Japan}, 
e-mail: {\tt ttro.skta.mki@gmail.com}
\end{small}
\end{flushleft}

\end{document}